\documentclass{ifacconf}

\usepackage{natbib}
\usepackage{amsmath}
\usepackage{amssymb}
\usepackage{mathtools}
\usepackage{array}
\usepackage{comment}
\usepackage{booktabs}
\usepackage{tikz}
\usetikzlibrary{arrows, arrows.meta}
\usepackage{algorithmic}
\usepackage{graphicx}
\usepackage{textcomp}
\usepackage{xcolor}
\usepackage{enumitem}
\usepackage{url}
\usepackage{siunitx}

\DeclareMathOperator{\Real}{Re} 

\newcommand{\suchthat}{\ifnum\currentgrouptype=16 \mathrel{}\middle|\mathrel{}\else\mid\fi}

\newcommand{\DominancyRefs}{\cite{Amrane2018Qualitative, BedouheneReal, Boussaada2018Dominancy, Boussaada2020Multiplicity, Boussaada2018Further, MazantiMultiplicity, Mazanti2020Qualitative, Mazanti2020Spectral, Boussaada2016Characterizing, Boussaada2016Tracking, Boussaada2016Multiplicity, Boussaada2018Towards, MBCBNC-Automatica-2022}}

\begin{document}

\begin{frontmatter}

\title{New Features of P$3\delta$ Software. \\ Insights and Demos\thanksref{footnoteinfo}}
\thanks[footnoteinfo]{The authors wish to acknowledge the work of the full P3$\delta$ development team, which, in addition to the authors, also include Yoann Audet, Thomas Charbonnet, Honor\'e Curlier,  Pierre-Henry Poret, Achrafy Said Mohamed and Franck Sim. The development of P3$\delta$ was also made possible thanks to the work of the Cyb'Air Association.}
\author[First,Second]{Islam Boussaada}
\author[First]{Guilherme Mazanti}
\author[First]{Silviu-Iulian Niculescu}
\author[Second]{Ayrton Hammoumou}
\author[Second]{Titouan Millet}
\author[Second]{Jayvir Raj}
\author[Second]{Julien Huynh}
\address[First]{Universit\'e Paris-Saclay, CNRS, CentraleSup\'elec, Inria, Laboratoire des Signaux et Syst\`emes, 91190, Gif-sur-Yvette, France. (e-mail: firstname.lastname@l2s.centralesupelec.fr).}
\address[Second]{IPSA, 63 Boulevard de Brandebourg, 94200 Ivry-sur-Seine, France}


\begin{abstract}
This paper presents the software entitled ``\emph{Partial Pole Placement via Delay Action},'' or ``P$3\delta$'' for short. P$3\delta$ is a Python software with a friendly user interface for the design of parametric stabilizing feedback laws with time-delays for dynamical systems. After recalling the theoretical foundation of the so-called ``\emph{Partial Pole Placement\/}'' methodology we propose as well the main features of the current version of P$3\delta$. We illustrate its use in feedback stabilization of several control systems operating under time delays.
\end{abstract}
\end{frontmatter}

\section{Introduction}

Time delays often occur in controlling dynamical systems, mainly due to the time required for acquiring, propagating, or processing information. It is commonly accepted that a delay in a control loop induces instability, oscillations and bad performance of the overall scheme. For instance, since the 30s, \cite{chp:36}, \cite{hpcs:37} showed the difficulty of handling delays in control loops for second- and third-order linear time-invariant (LTI) dynamical systems and one of the natural ideas was to compensate it by pre-correction (see, e.g., \cite{porter:52}).  

However, as noted and briefly discussed in \cite{Sipahi2011Stability}, in some cases, the delay can have a stabilizing effect. Furthermore, it has been emphasized in \cite{Suh1979Proportional} and \cite{Atay1999Balancing} that one may replace the classical proportional-derivative controller by a proportional-delayed controller, using the so-called ``average derivative action'' due to the delay. Finally, regarding the beneficial effect of the delay, closed-loop stability may be guaranteed for some dynamical systems subject to input delays precisely by the existence of such delays, as pointed out in \cite{Niculescu2010Delay} and the references therein in controlling oscillators by delayed output feedback. 

The above reasons explain, in part, an abundant literature on such topics, such as \cite{Gu2003Stability, Michiels2014Stability, Stepan1989Retarded} and the references therein. It should be mentioned that time-delay dynamical systems are infinite-dimensional systems and there exist several ways to represent their dynamics. In the sequel,
their dynamics are represented by delay-differential equations (DDEs). For an introduction on DDEs, we refer to \cite{Hale1993Introduction}.

The software\footnote{An acronym for \emph{Partial Pole Placement via Delay Action\/}} P3$\delta$ was introduced in \cite{Boussaada-ICSTCC2020P3D,Boussaada-TDS2021}. The main intention of the authors is to help users interested in stability analysis and stabilization of dynamical LTI systems in the presence of a single delay in closed loop.  P3$\delta$ makes use of the so-called \emph{mul\-ti\-pli\-city-in\-duced-do\-mi\-nancy (MID)\/} and \emph{partial pole placement\/} methods. More precisely, MID is an unexpected spectral property stating that, in some cases, the characteristic root with maximal multiplicity defines the spectral abscissa of the corresponding characteristic function, i.e., the rightmost root of the spectrum. In control, MID opened a new and interesting perspective relying on the idea of a \emph{partial pole placement\/} by reinforcing the use of the delay as a \emph{control parameter}. For a deeper discussion on these methods, we refer to \DominancyRefs.

The software P3$\delta$ covers DDEs of retarded or neutral type\footnote{For the classification of DDEs, we refer to \cite{Hale1993Introduction}.} with a single time delay, under the form
\begin{multline}
\label{MainEqn}
\textstyle y^{(n)}(t) + a_{n-1} y^{(n-1)}(t) + \dotsb + a_0 y(t) \\
\textstyle + b_m y^{(m)}(t - \tau) + \dotsb + b_0 y(t - \tau) = 0,
\end{multline}
under appropriate initial conditions, where $\tau > 0$ is the positive delay, $y$ is the real-valued unknown function, $n$ and $m$ are nonnegative integers with $n \geq m$, and $a_0, \dotsc, a_{n-1}, b_0, \dotsc, b_m$ are real coefficients. 

To address the stability analysis of LTI DDEs, the software P3$\delta$ relies on spectral methods (see, e.g., \cite{Hale1993Introduction, Michiels2014Stability}), which consist on the study of the complex roots of a \emph{characteristic function} of the system. The characteristic function $\Delta : \mathbb{C} \to \mathbb{C}$ of \eqref{MainEqn} is given by
\begin{equation}
\label{Delta}
\Delta(s) = s^n + \sum_{k=0}^{n-1} a_k s^k + e^{-s \tau} \sum_{k=0}^m b_k s^k,
\end{equation}
and \eqref{MainEqn} is exponentially stable if and only if the \emph{spectral abscissa} $\gamma = \sup\{\Real s \suchthat \Delta(s) = 0\}$ satisfies $\gamma < 0$.

P3$\delta$ considers that $b_0, \dotsc, b_m$ are free parameters. In its ``Generic'' mode, $a_0, \dotsc, a_{n-1}$ are also assumed to be free and $\tau$ is fixed while, in its ``Control-oriented'' mode, $a_0, \dotsc, a_{n-1}$ are assumed fixed and $\tau$ can be assumed either free or fixed. The ``Control-oriented'' mode is usually suitable for control applications, as illustrated in Section~\ref{SecP3Delta}.

The strategy used by P3$\delta$ to stabilize a time-delay system is to tune its free parameters to assign finitely many roots while also guaranteeing that the rightmost root on the complex plane is among the chosen ones. Two main properties define such a strategy: (i) first, assigning a real root of maximal multiplicity and proving that this root is necessarily the rightmost root of the characteristic quasipolynomial, a property which has been named \emph{mul\-ti\-pli\-city-in\-duced-do\-mi\-nancy}, or MID for short, and (ii) second, assigning a certain amount of real roots, typically equally spaced for simplicity, and proving that the rightmost root among them is also the rightmost root of the characteristic quasipolynomial, a property which has been named \emph{coexisting real roots-in\-duced-do\-mi\-nancy}, or CRRID for short, in \cite{Amrane2018Qualitative,BedouheneReal}.

The MID property for \eqref{MainEqn} was shown, for instance, in \cite{Boussaada2018Further} in the case $(n, m) = (2, 0)$, in \cite{Boussaada2020Multiplicity} in the case $(n, m) = (2, 1)$ (see also \cite{Boussaada2018Dominancy}), in \cite{MazantiMultiplicity} in the case of any positive integer $n$ and $m = n-1$ (see also \cite{Mazanti2020Qualitative}), and more recently in \cite{BMN-CRAS-2022} for arbitrary $n\geq m$. It was also studied for neutral systems of orders $1$ and $2$ in \cite{MBCBNC-Automatica-2022,Benarab2020MID} and extended to complex conjugate roots of maximal multiplicity in \cite{Mazanti2020Spectral}. The CRRID property was shown, for instance, in \cite{Amrane2018Qualitative} in the cases $(n, m) = (2, 0)$ and $(n, m) = (1, 0)$, and in \cite{BedouheneReal} in the case of any positive integer $n$ and $m = 0$.

In all the above cases, the maximal multiplicity of a real root or, equivalently, the maximal number of coexisting simple real roots is $n + m + 1$. The idea to exploit the nature of (real or complex) open-loop roots in control design was proposed for second-order systems in \cite{Boussaada2020Multiplicity} and extended for arbitrary order systems with real-rooted plants in \cite{Balogh20,Balogh21}.

As mentioned earlier, P3$\delta$ allows for the parametric design of stabilizing feedback laws with time-delays by exploiting the MID and CRRID properties briefly presented above.
The present paper describes the new functionalities of P$3\delta$ and  provides some illustrative examples for its use.

\section{New features of P3$\delta$}

\subsection{New features of the online version}

\begin{figure*}[ht]
\centering
\begin{minipage}[b]{0.333\textwidth}
\centering

\includegraphics[width=0.84\textwidth]{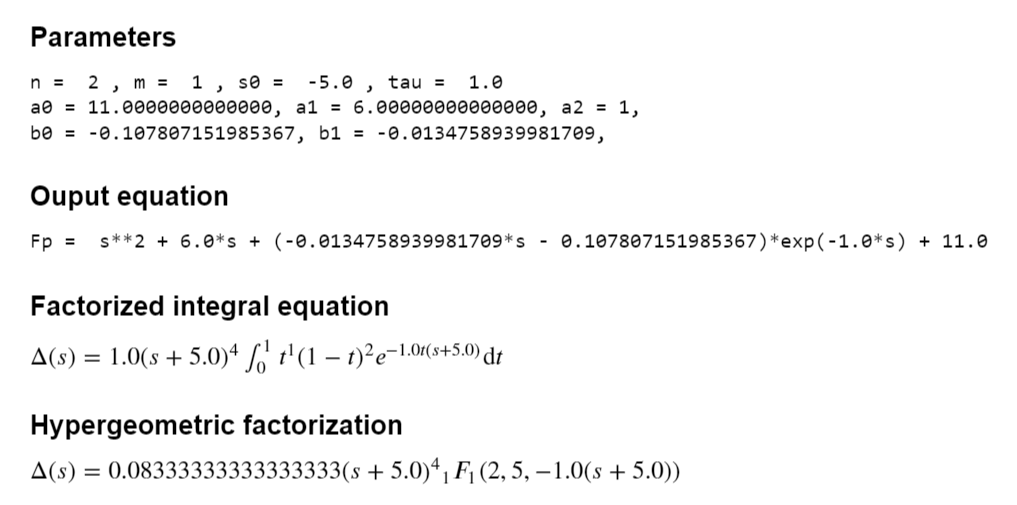}

(a)

\end{minipage}\begin{minipage}[b]{0.333\textwidth}
\centering

\includegraphics[width=0.84\textwidth]{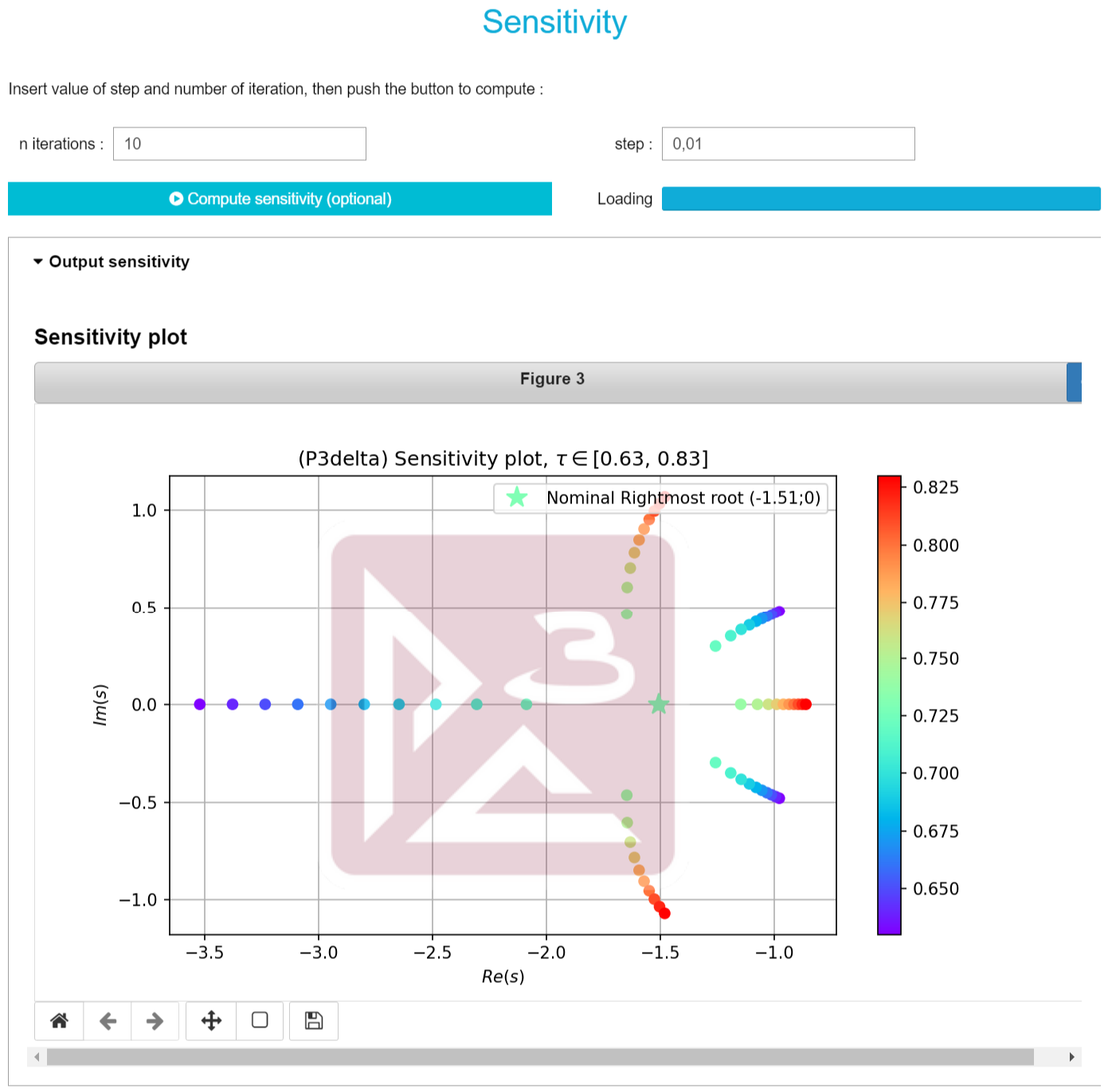}

(b)

\end{minipage}\begin{minipage}[b]{0.333\textwidth}
\centering

\includegraphics[width=0.84\textwidth]{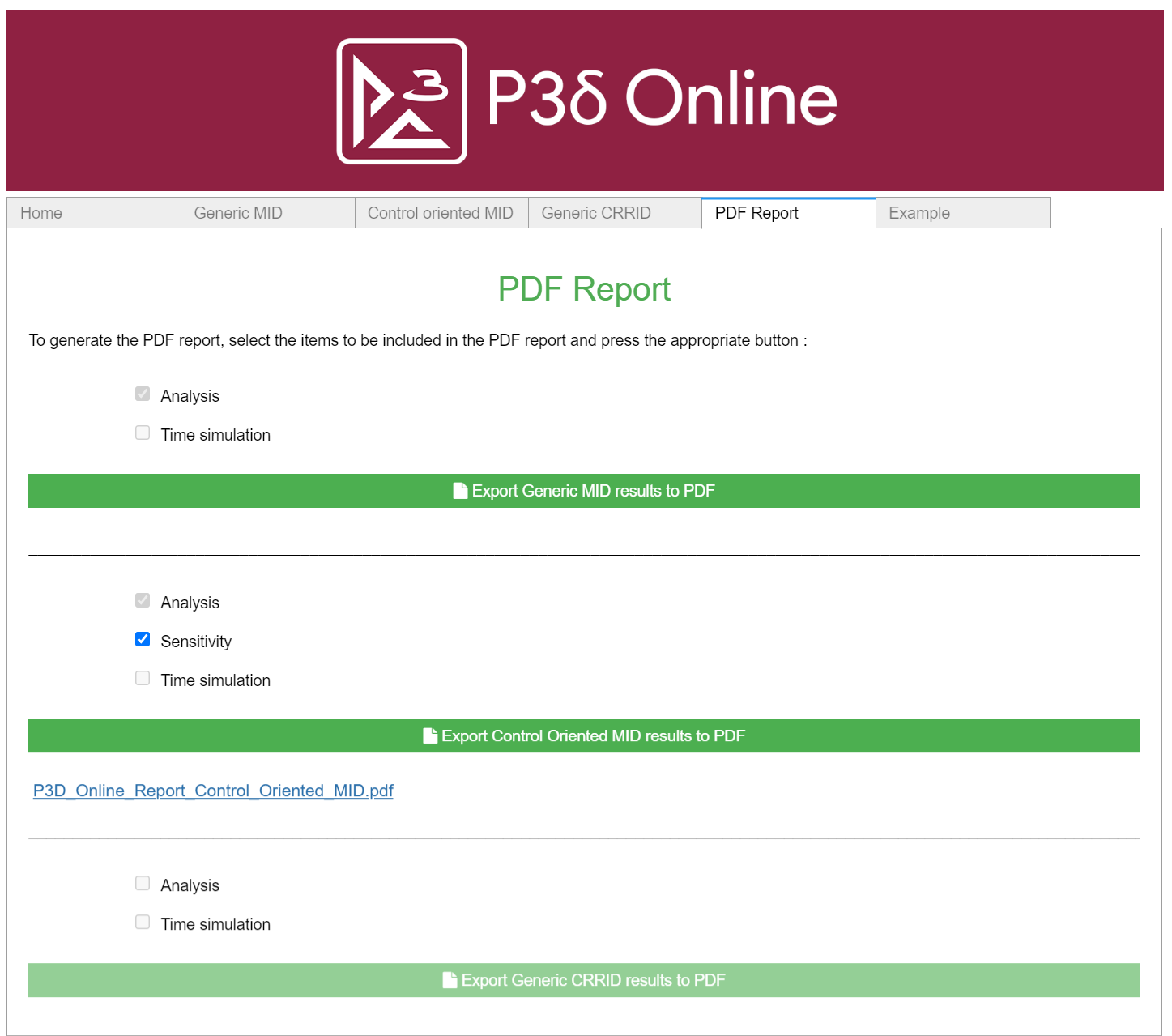}

(c)
\end{minipage}
\caption{Illustration of the different new features of P3$\delta$ Online.}
\label{FigNewFeatures}
\end{figure*}

After a first iteration of the software, P3$\delta$ Online was completed with new features (Figure~\ref{FigNewFeatures}) that enrich the software and improve user experience. The one-click online version of P3$\delta$ is still hosted on \emph{Binder}, which provides a personalized computing environment directly from a \emph{GitHub} repository. The P3$\delta$ team continued to develop the online software in Python, using the dynamic \emph{Jupyter Notebook} format and with an user interface built using interactive widgets from Python's \texttt{ipywidgets} module. P3$\delta$ is freely available for download on \url{https://cutt.ly/p3delta}, where installation instructions, video demonstrations, and the user guide are also available.\footnote{Interested readers may also contact directly any of the authors of the paper.}

P3$\delta$ Online is based on the program of the executable version and enriched with exclusive features  to this version. The online software includes features from the ``Generic MID'', ``Control-oriented MID'', and ``Generic CRRID'' modes of P3$\delta$ described in \cite{Boussaada-TDS2021, Boussaada-ICSTCC2020P3D}.

In the ``Generic MID'' mode, the online version of the software returns the spectrum distribution as well as a normalized quasipolynomial which admits a root of multiplicity $n+m+1$ at the origin\footnote{In this case, the multiplicity coincides with the degree of the quasipolynomial $\Delta$, and represents the maximal allowable multiplicity.}. In the ``Control-oriented MID'' mode, the online version of software returns the admissibility region, a normalized quasipolynomial which admits a root of multiplicity $m+2$ at the origin, and an illustration of the bifurcation of the root of multiplicity $m+2$ with respect to variations of the delay $\tau$.

The first new feature of P3$\delta$ Online is the addition of a ``Home'' tab that contains information about the project and the software settings. Two settings are currently available: the appearance of notifications and the activation of the software limits.

In the ``Generic MID'' mode, the spectral distribution analysis now provides two new forms of the output equation, Factorized integral equation and Hypergeometric factorization. In the ``Control-oriented MID'' mode, it is now possible to choose the value of the discretization step as well as the number of iterations for the $\tau$ sensitivity plot.

Finally, the P3$\delta$ team has developed a new feature to export the results obtained in the software in the format of a report automatically written in a PDF file. The user has the choice of which calculation mode results they want to export and can save the report directly on their computer.

\subsection{Assignment admissibility region}
\label{AdmissibilityRegi}

Given a system, the definition of valid hyper-parameters ($s_0$ and $\tau$) can sometimes be difficult. In order to help users getting a better idea of the admissibility ranges of $s_0$ and $\tau$, P3$\delta$ has an admissibility region plotting feature which will be described in this section. Note that this feature is only available when using ``Control-oriented MID'' as it currently is the only mode with constraints on $(s_0, \tau)$.

More precisely, given $a_0, \dotsc, a_{n-1}$, the admissibility region is defined as the set of pairs $(s_0, \tau) \in \mathbb R \times (0, +\infty)$ for which there exist real coefficients $b_0, \dotsc, b_m$ such that $s_0$ is a root of $\Delta$ of multiplicity at least $m + 2$ when the delay is $\tau$. To compute such a region, we first express the coefficients $(b_i)_{0 \leq i \leq m}$ in terms of $(a_i)_{0\leq i \leq n-1}$, $s_0$, and $\tau$ by using the $m+1$ equations $\Delta^{(k)}(s_0) = 0$, $k \in \{0, \dotsc, m\}$. As those equations are linear in the $m+1$ variables $b_0, \dotsc, b_m$, this linear system admits a unique solution. We then replace these expressions of $(b_i)_{0 \leq i \leq m}$ in the equation $\Delta^{(m+1)}(s_0) = 0$, obtaining an algebraic relation between $s_0$, $\tau$, and the coefficients $(a_i)_{0 \leq i \leq n-1}$. By construction, this is a necessary and sufficient condition for $s_0$ to be a root of multiplicity at least $m+2$ of $\Delta$ and, since the coefficients $(a_i)_{0 \leq i \leq n-1}$ are known in this mode, this algebraic equation characterizes the admissibility region.

In P3$\delta$, the computations leading to this admissibility region is done by following the previously mentioned steps in a symbolic way using the \texttt{sympy} package. Only the part of the admissibility region in the rectangle $[s_{0, \min}, 0] \times [0, \tau_{\max}]$ is displayed, where $s_{0, \min} < 0$ and $\tau_{\max} > 0$ are values selected by the user.

\section{Applications of P3$\delta$}
\label{SecP3Delta}

To illustrate P3$\delta$, we present some applications of the software.

\subsection{Harmonic oscillator}
\label{SecGenericMID}

Consider a controlled harmonic oscillator described by 
\begin{equation}
\label{eq:oscillator}
y''(t) + y(t) = u(t),
\end{equation}
where $y(t) \in \mathbb R$ is the instantaneous state of the oscillator available to measurement, the control input $u(t)$ corresponds to the applied force, and the coefficients of the equation were normalized. We assume that the control input is given by a delayed proportional-derivative controller
\begin{equation}
\label{eq:oscillator-feedback}
u(t) = -\beta y(t - \tau) - \alpha y'(t - \tau),
\end{equation}
where $\alpha$ and $\beta$ are the coefficients of the controller and $\tau > 0$ is the delay. The characteristic equation of the closed-loop system is thus
\begin{equation}
\label{eq:delta-harmonic}
\Delta (s) = s^2 + 1 + (\beta +\alpha s) e^{-s\tau}.
\end{equation}
The polynomial corresponding to the non-delayed term in \eqref{eq:delta-harmonic} is of degree $n = 2$, while that corresponding to the delayed term is of degree $m = 1$. The degree of the quasipolynomial $\Delta$ is thus $n + m + 1 = 4$.

\begin{figure}[ht]
\centering
\includegraphics[width=0.84\columnwidth]{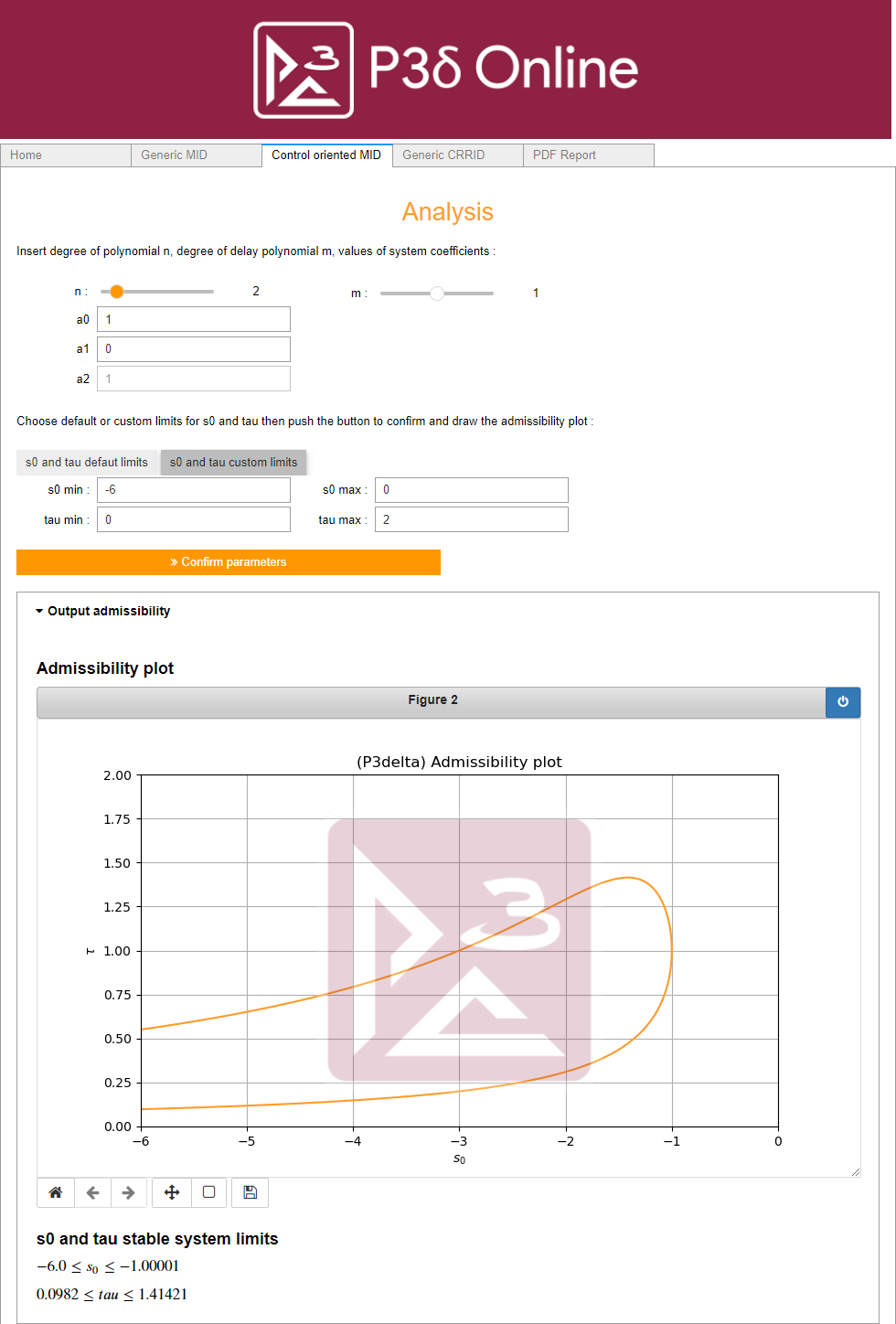}
\caption{Data input and admissibility plot for the harmonic oscillator.}
\label{fig:oscillator-1}
\end{figure}

Let us use P3$\delta$ in order to place a root of multiplicity $m+2 = 3$ at some $s_0 \in \mathbb R$. One should first input into P3$\delta$ the values of the known coefficients of the qua\-si\-po\-ly\-no\-mial $\Delta$ from \eqref{eq:delta-harmonic}, i.e., the coefficients of the polynomial corresponding to the non-delayed term. After introducing the data corresponding to \eqref{eq:delta-harmonic} into P3$\delta$, one obtains the plot of the admissibility region, as in Figure~\ref{fig:oscillator-1}.

\begin{figure}[ht]
\centering
\includegraphics[width=0.84\columnwidth]{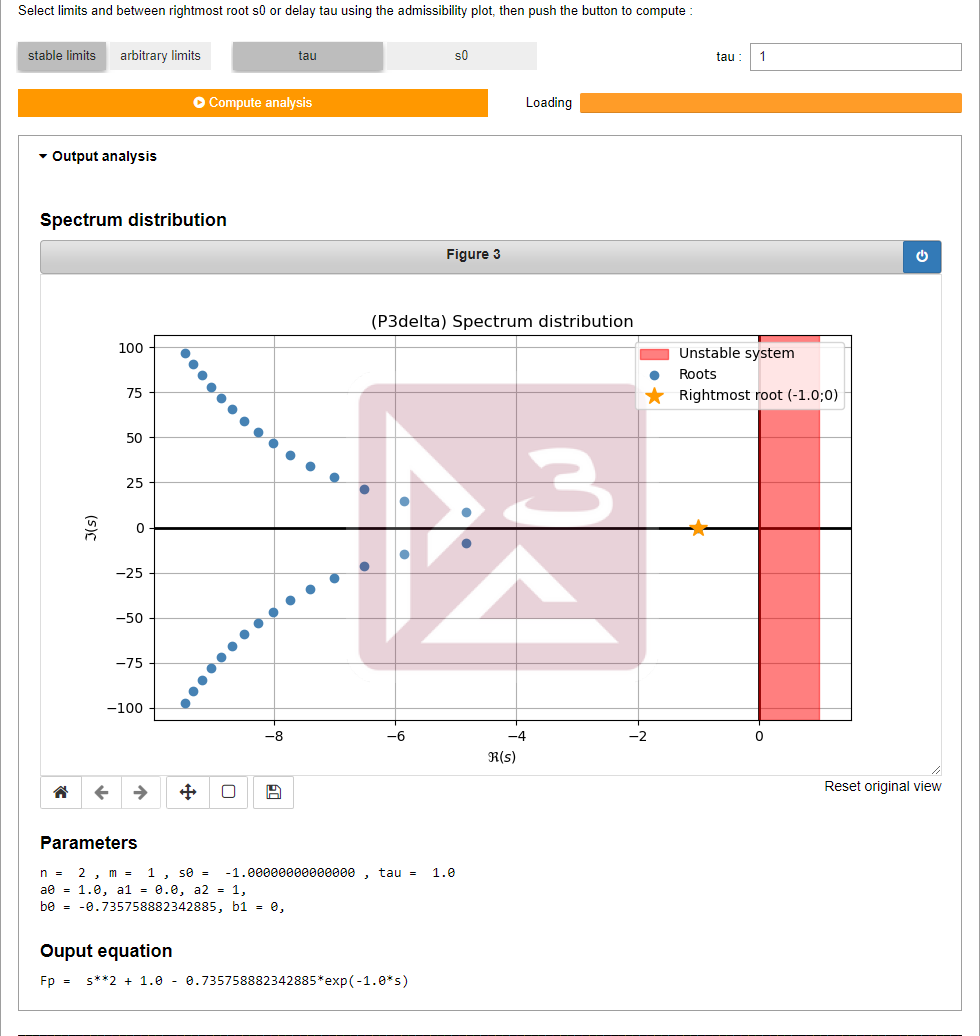}
\caption{Computed values and spectrum distribution for the harmonic oscillator.}
\label{fig:oscillator-2}
\end{figure}

Based on Figure~\ref{fig:oscillator-1}, one observes that the method used by P3$\delta$ allows for the stabilization of the system if the delay is at most approximately $1.4$. After inputting the value of the delay, P3$\delta$ will compute the root $s_0$ and the coefficients $\alpha$ and $\beta$, and trace the spectrum distribution of \eqref{eq:delta-harmonic} with the derived values. In the present example, the corresponding screen of P3$\delta$ is shown in Figure~\ref{fig:oscillator-2}, where, for $\tau = 1$, one obtains $s_0 = -1$, $\alpha = 0$, and $\beta \approx -0.7358$.

\subsection{Inverted pendulum}
\label{Secinvertpendulum}

As a second example of application of P3$\delta$, we consider the stabilization of an inverted pendulum from Figure~\ref{fig:pendulum}, in which a stick of mass $m$ and length $\ell$ is placed over a cart and can rotate freely around the attachment point $A$. The cart can move along a rail and an external force $F$ is applied on the cart, and we assume that the mass of the cart is negligible compared to the mass of the stick. We denote by $\varphi$ the angle between the stick and the upward vertical direction. In this case, following \cite{Molnar2021Calculation}, the linearization of the dynamics of $\varphi$ around the unstable equilibrium $\varphi = 0$ is
\begin{equation}
\label{eq:pendulum}
\textstyle \varphi''(t) - \frac{m g \ell}{2I}\varphi(t) = \frac{\ell}{2I}F(t),
\end{equation}
where $I = \frac{1}{12} m \ell^2$ is the moment of inertia of the stick and $g$ is the gravitational acceleration. For the numerical application, we will consider 
$m = 10$~kg
and
$\ell = 10$~m,
in which case
$I \approx 83.33$~kg~m\textsuperscript{2}.

\newcommand{\pendangle}{15}
\begin{figure}[ht]
\centering
\resizebox{0.65\columnwidth}{!}{
\begin{tikzpicture}
\draw (-3.5, 0) -- (3.5, 0);
\draw (-3.5, -0.5) -- (-3.5, 0.5);
\draw (3.5, -0.5) -- (3.5, 0.5);
\foreach \pos in {-0.5, -0.4, ..., 0.6} {
    \draw (-3.55, {\pos-0.05}) -- (-3.5, \pos);
    \draw (3.5, \pos) -- (3.55, {\pos+0.05});
}
\draw[dashed] (0, 0) -- (0, 4);
\draw[fill=white] (-1, -0.5) rectangle (1, 0.5);
\draw[fill=white] (0, 0) -- ++({180-\pendangle}:0.1) -- ++({90-\pendangle}:4) -- ++(-\pendangle:0.2) -- node[midway, right] {$m$, $\ell$} ++({-90-\pendangle}:4) -- cycle;
\draw[fill=white] (0, 0) circle[radius=0.1];
\draw[->] (0, 3) arc [start angle=90, end angle={90-\pendangle}, radius=3] node[midway, above] {$\varphi$}; 
\draw (0, 0) node[below right] {$A$};
\draw[->, blue] (-2, -0.25) node[below]{$F$} -- (-0.2, -0.25);
\end{tikzpicture}
}
\caption{Inverted pendulum on a cart.}
\label{fig:pendulum}
\end{figure}
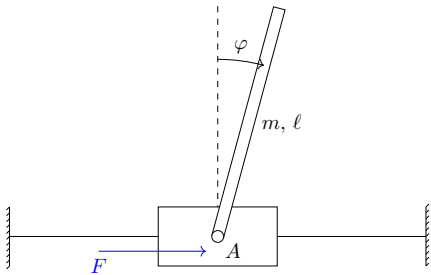

We wish to stabilize the origin of \eqref{eq:pendulum} by a delayed PD feedback of the angular position, i.e., we wish to apply
\begin{equation*}
    F(t) = -K_p\varphi(t-\tau) - K_d \varphi'(t-\tau),
\end{equation*}
yielding a closed-loop system with characteristic equation
\begin{equation}
\label{eq:Delta-pendulum}
\textstyle \Delta(s) = s^2 - \frac{m g \ell}{2I} + \left(\frac{K_p \ell}{2I} + \frac{K_d \ell}{2I} s \right) e^{-s\tau}.
\end{equation}
Note that, here, the delay $\tau$ can be seen as a design parameter, together with $K_p$ and $K_d$.

\begin{figure}
    \centering
    \includegraphics[width=0.84\columnwidth]{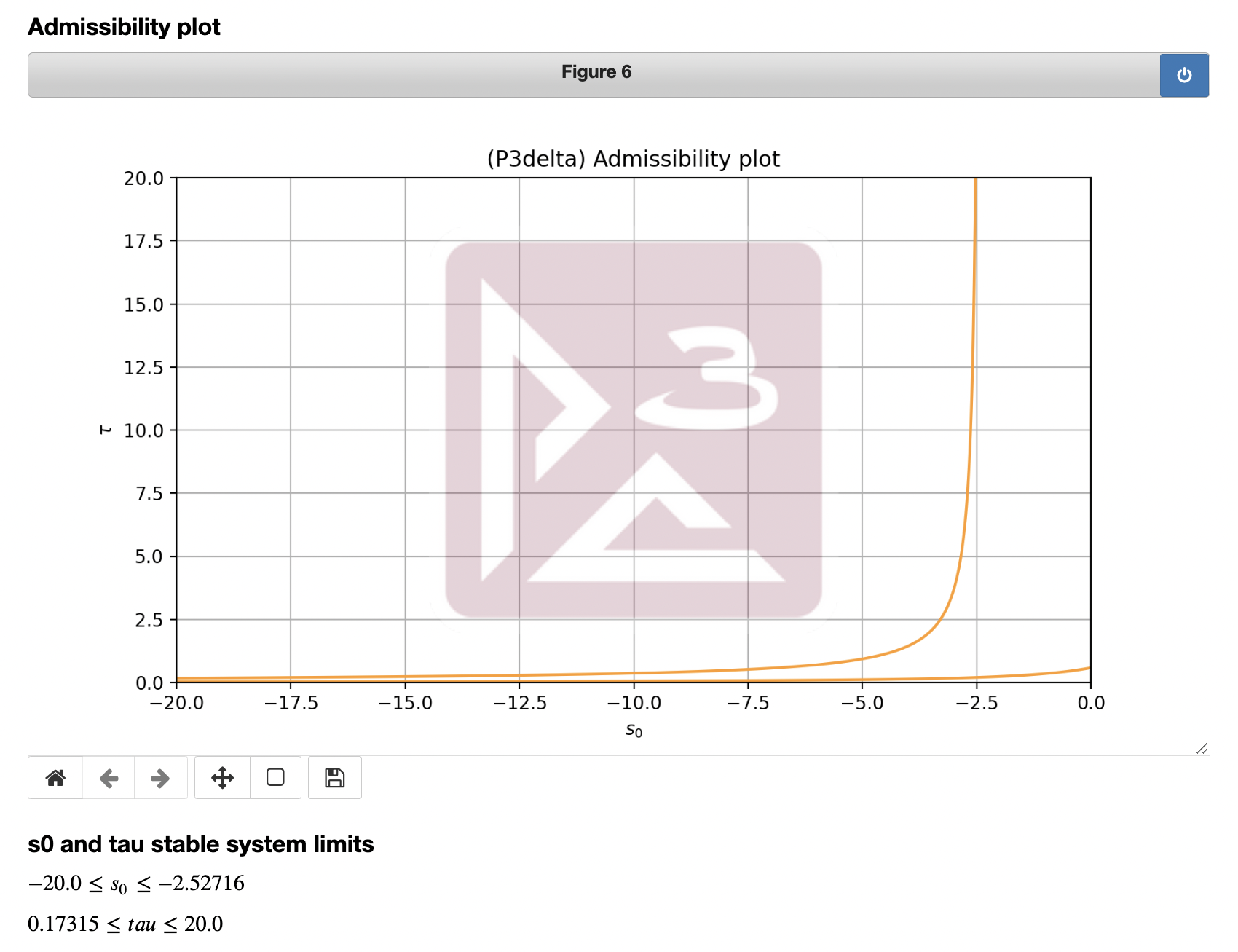}
    \caption{Admissibility plot for the inverted pendulum.}
    \label{fig:pendulum-admissibility}
\end{figure}

\begin{figure}
    \centering
    \includegraphics[width=0.84\columnwidth]{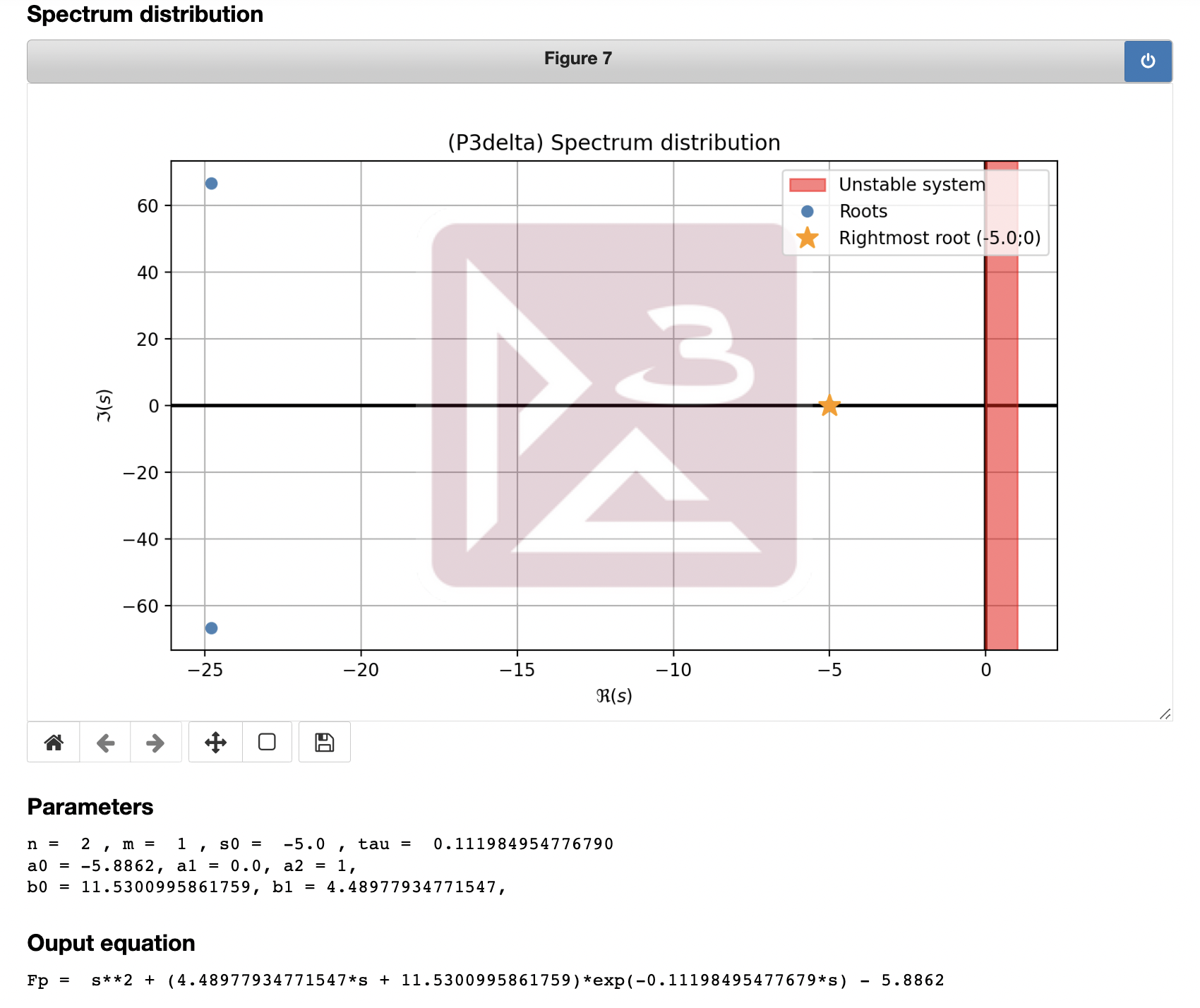}
    \caption{Computed values and spectrum distribution for the inverted pendulum.}
    \label{fig:pendulum-spectrum}
\end{figure}

By introducing the data corresponding to \eqref{eq:Delta-pendulum} into P3$\delta$ ``Control-oriented MID'' mode, we obtain the admissibility plot from Figure~\ref{fig:pendulum-admissibility}. Based on this plot, we decide to assign a root of multiplicity $3$ of \eqref{eq:Delta-pendulum} at $s_0 = -5$. P3$\delta$ then gives the result shown in Figure~\ref{fig:pendulum-spectrum}. From the computed values of the coefficients shown in Figure~\ref{fig:pendulum-spectrum}, one deduces that the delay corresponding to $s_0 = -5$ is $\tau \approx 0.112$ and the coefficients of \eqref{eq:Delta-pendulum} should be $\frac{K_p \ell}{2I} \approx 11.53$ and $\frac{K_d \ell}{2I} \approx 4.4898$, yielding the parameters $K_p \approx 192.16$ and $K_d \approx 74.83$.

\subsection{Transonic flow in a wind tunnel}
\label{Transflow}

Consider now a wind tunnel in which a cold fluid is set to motion at high speed. The control of the velocity of such a fluid around an equilibrium state can be described by the system
\begin{equation}
\label{eq:flow}
\left\{
\begin{aligned}
\kappa m'(t) + m(t) & = k \theta(t - \tau_0), \\
\theta''(t) + 2 \zeta \omega \theta'(t) + \omega^2 \theta(t) & = u(t),
\end{aligned}
\right.
\end{equation}
where $m$ denotes the deviation of the Mach number of the fluid with respect to its equilibrium state, $\kappa$ and $k$ are constants depending on the characteristics of the fluid and the desired equilibrium state, $\theta$ is the angle of a guide vane driving the velocity of the fluid, $\tau_0 > 0$ is a delay depending only on the temperature of the fluid, $\zeta$ and $\omega$ are parameters of the dynamics of the guide vane angle, and $u$ is a control input. The above model comes from \cite{Armstrong1981Application} and its stabilization was previously discussed in \cite{MazantiMultiplicity} under the assumption that one may choose $\zeta$ and $\omega$. We consider here a more realistic situation in which $\zeta$ and $\omega$ are fixed, and a feedback law of the form 
\[
u(t) = -\beta m(t - \tau_1) - \alpha_0 \theta(t - \tau_0 - \tau_1) - \alpha_1 \theta'(t - \tau_0 - \tau_1),
\]
where $\tau_1 \geq 0$ is a new delay, which can be seen as a design parameter and should be at least equal to the delay for measuring $m$. For simplicity, we denote $\tau = \tau_0 + \tau_1$.  Inserting this control law into \eqref{eq:flow} and performing straightforward algebraic manipulations, one deduces that $\theta$ verifies the third-order differential equation
\begin{multline*}
\textstyle \theta'''(t) + \left(2 \zeta \omega + \frac{1}{\kappa}\right) \theta''(t) + \left(\omega^2 + \frac{2 \zeta \omega}{\kappa}\right) \theta'(t) + \frac{\omega^2}{\kappa} \theta(t) \\
\textstyle {} + \alpha_1 \theta''(t - \tau) + \left(\alpha_0 + \frac{\alpha_1}{\kappa}\right) \theta'(t - \tau) + \frac{\alpha_0 + \beta k}{\kappa} \theta(t - \tau) = 0,
\end{multline*}
and hence the closed-loop characteristic quasipolynomial $\Delta$ of \eqref{eq:flow} is
\begin{equation}
\label{eq:delta-flow}
\begin{split}
\textstyle \Delta(s) = {} & \textstyle s^3 + \left(2 \zeta \omega + \frac{1}{\kappa}\right) s^2 + \left(\omega^2 + \frac{2 \zeta \omega}{\kappa}\right) s + \frac{\omega^2}{\kappa} \\
& \textstyle {} + (\gamma_2 s^2 + \gamma_1 s + \gamma_0)e^{-s\sigma},
\end{split}
\end{equation}
where $\gamma_2 = \alpha_1$, $\gamma_1 = \alpha_0 + \frac{\alpha_1}{\kappa}$, and $\gamma_0 = \frac{\alpha_0 + \beta k}{\kappa}$.

\begin{figure}[ht]
    \centering
    \includegraphics[width=0.84\columnwidth]{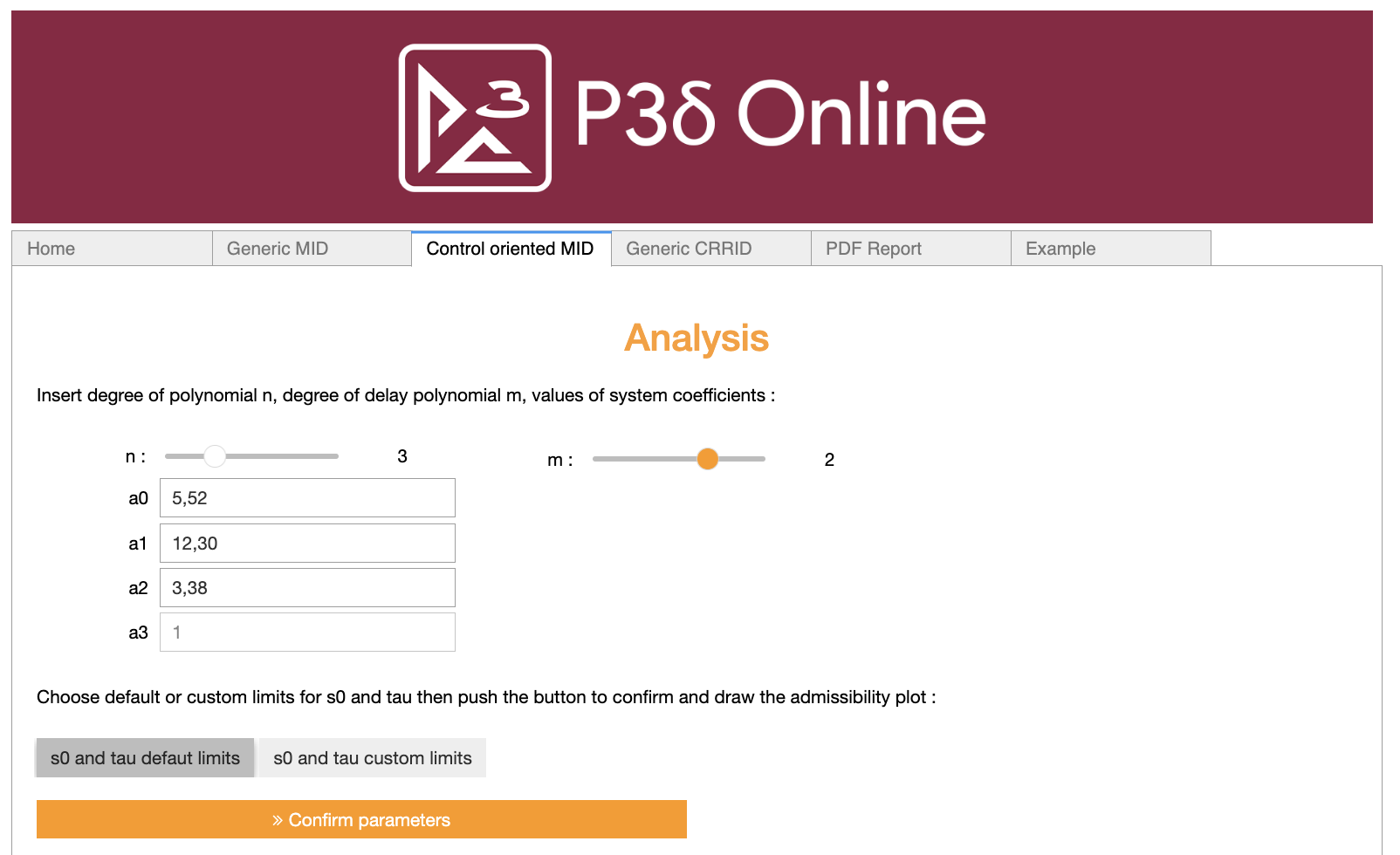}\\
    \caption{Data input for the transonic flow model.}
    \label{fig:flow-data}
\end{figure}

\begin{figure}[ht]
    \centering
    \includegraphics[width=0.84\columnwidth]{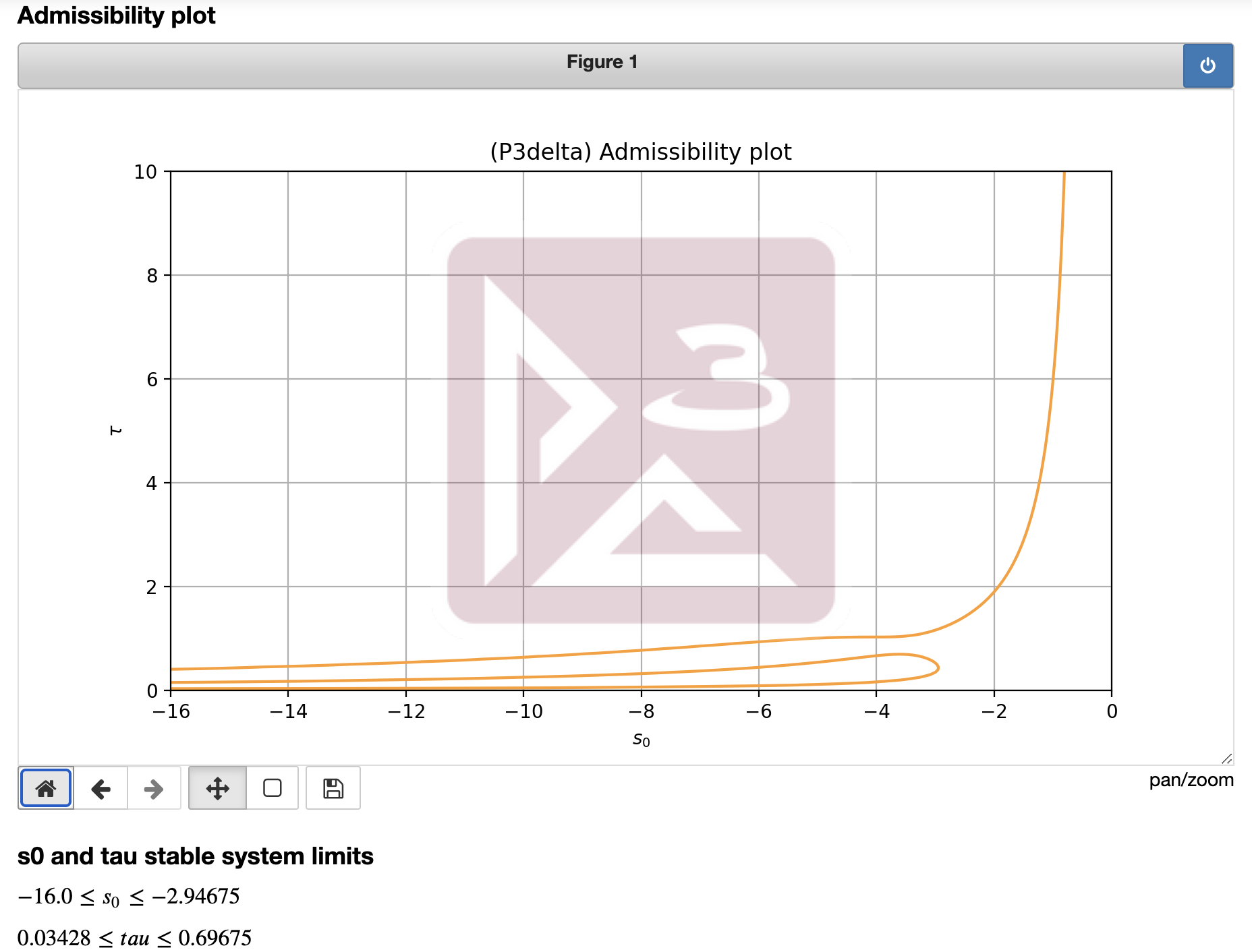}\\
    \caption{Admissibility plot for the transonic flow model.}
    \label{fig:flow-admissible}
\end{figure}

As a numerical application, we consider the linearization around the steady state with Mach number $0.84$ and air temperature 
$166$~K,
in which case, as reported in \cite{Armstrong1981Application}, the system parameters are
$\kappa = 1.964$~s,
$k = -0.67036$~rad$^{-1}$,
$\tau_0 = 0.33$~s,
$\zeta = 0.4368$, and
$\omega = 3.292$~rad~s$^{-1}$.
We can hence insert the known coefficients of \eqref{eq:delta-flow} into P3$\delta$, as shown in Figure~\ref{fig:flow-data}, and obtain the admissibility plot from Figure~\ref{fig:flow-admissible}.
    
\begin{figure}[ht]
    \centering
    \includegraphics[width=0.84\columnwidth]{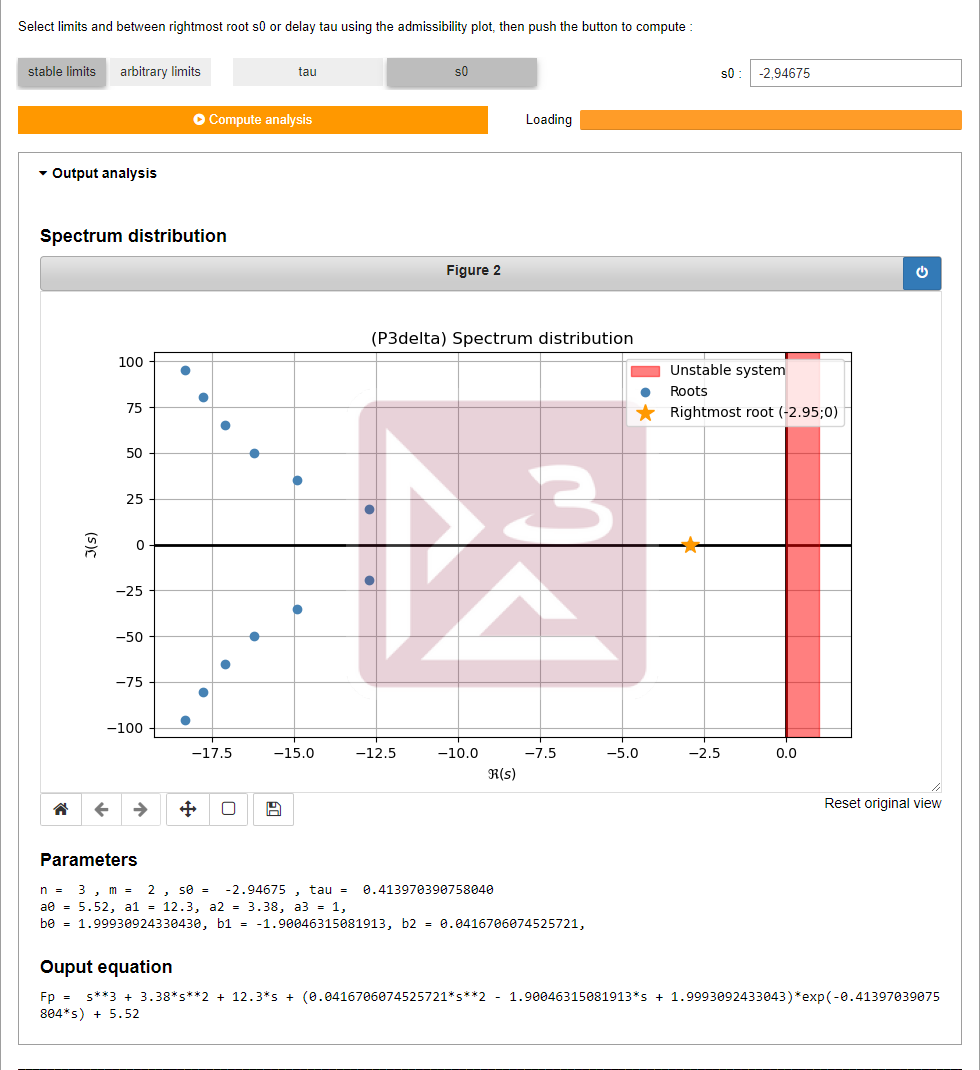}\\
    \caption{Computed values and spectrum distribution for the transonic flow model.}
    \label{fig:flow-spectrum}
\end{figure}

Based on Figure~\ref{fig:flow-admissible}, we choose in this example $s_0 = -2.94675$ as a root of multiplicity $5$ of the quasipolynomial $\Delta$, and P3$\delta$ provides the output shown in Figure~\ref{fig:flow-spectrum}. In particular, we deduce the numerical value
$\tau \approx 0.4140$~s,
yielding
$\tau_1 \approx 0.0840$~s,
as well as $\gamma_2 \approx 0.04167$, $\gamma_1 \approx -1.9005$, and $\gamma_0 \approx 1.9993$, yielding $\beta \approx -2.9908$, $\alpha_0 \approx 1.9217$, and $\alpha_1 \approx 0.04167$.

\section{Concluding remarks and future work}

By exploiting the properties of the MID and the CRRID, P3$\delta$ allows for the design of feedback control laws for real applications. The main novelty of the version discussed in the paper is the improvement of the graphic interface for the online version of the software. Inspired by the ``Control-oriented MID'' mode, the team is working on a ``Control-oriented CRRID'' mode. In addition, an executable version for macOS is currently under development.

\section*{Acknowledgments}
This work is partially supported by a grant by the French National Research Agency (ANR) as part of the ``Investissement d'Avenir'' program, through the iCODE project funded by the IDEX Paris-Saclay, ANR-11-IDEX0003-02. The authors also acknowledge the support of Institut Polytechnique des Sciences Avanc\'ees (IPSA).

\bibliography{p3delta}

\begin{thebibliography}{32}
\providecommand{\natexlab}[1]{#1}
\providecommand{\url}[1]{\texttt{#1}}
\providecommand{\urlprefix}{URL }
\expandafter\ifx\csname urlstyle\endcsname\relax
  \providecommand{\doi}[1]{doi:\discretionary{}{}{}#1}\else
  \providecommand{\doi}{doi:\discretionary{}{}{}\begingroup
  \urlstyle{rm}\Url}\fi

\bibitem[{Amrane et~al.(2018)Amrane, Bedouhene, Boussaada, and
  Niculescu}]{Amrane2018Qualitative}
Amrane, S., Bedouhene, F., Boussaada, I., and Niculescu, S.I. (2018).
\newblock On qualitative properties of low-degree quasipolynomials: further
  remarks on the spectral abscissa and rightmost-roots assignment.
\newblock \emph{Bull. Math. Soc. Sci. Math. Roumanie (N.S.)}, 61(109)(4),
  361--381.

\bibitem[{Armstrong and Tripp(1981)}]{Armstrong1981Application}
Armstrong, E.S. and Tripp, J.S. (1981).
\newblock An application of multivariable design techniques to the control of
  the {N}ational {T}ransonic {F}acility.
\newblock Technical Report 1887, NASA.

\bibitem[{Atay(1999)}]{Atay1999Balancing}
Atay, F.M. (1999).
\newblock Balancing the inverted pendulum using position feedback.
\newblock \emph{Appl. Math. Lett.}, 12(5), 51--56.

\bibitem[{Balogh et~al.(2020)Balogh, Insperger, Boussaada, and
  Niculescu}]{Balogh20}
Balogh, T., Insperger, T., Boussaada, I., and Niculescu, S.I. (2020).
\newblock Towards an {MID}-based delayed design for arbitrary-order dynamical
  systems with a mechanical application.
\newblock \emph{IFAC-PapersOnLine}, 53(2), 4375--4380.
\newblock 21th IFAC World Congress.

\bibitem[{Balogh et~al.(2022)Balogh, Boussaada, Insperger, and
  Niculescu}]{Balogh21}
Balogh, T., Boussaada, I., Insperger, T., and Niculescu, S.I. (2022).
\newblock Conditions for stabilizability of time-delay systems with real-rooted
  plant.
\newblock \emph{Internat. J. Robust Nonlinear Control}, 32(6), 3206--3224.

\bibitem[{Bedouhene et~al.(2020)Bedouhene, Boussaada, and
  Niculescu}]{BedouheneReal}
Bedouhene, F., Boussaada, I., and Niculescu, S.I. (2020).
\newblock Real spectral values coexistence and their effect on the stability of
  time-delay systems: Vandermonde matrices and exponential decay.
\newblock \emph{Comptes Rendus. Math\'ematique}, 358(9-10), 1011--1032.

\bibitem[{Benarab et~al.(2020)Benarab, Boussaada, Trabelsi, Mazanti, and
  Bonnet}]{Benarab2020MID}
Benarab, A., Boussaada, I., Trabelsi, K., Mazanti, G., and Bonnet, C. (2020).
\newblock The {MID} property for a second-order neutral time-delay differential
  equation.
\newblock In \emph{2020 24th International Conference on System Theory, Control
  and Computing (ICSTCC)}, 202--207.

\bibitem[{Boussaada et~al.(2022)Boussaada, Mazanti, and
  Niculescu}]{BMN-CRAS-2022}
Boussaada, I., Mazanti, G., and Niculescu, S.I. (2022).
\newblock The generic multiplicity-induced-dominancy property from retarded to
  neutral delay-differential equations: {When} delay-systems characteristics
  meet the zeros of {Kummer} functions.
\newblock \emph{Comptes Rendus. Math.}, 360, 349--369.

\bibitem[{Boussaada et~al.(2020{\natexlab{a}})Boussaada, Mazanti, Niculescu,
  Huynh, Sim, and Thomas}]{Boussaada-ICSTCC2020P3D}
Boussaada, I., Mazanti, G., Niculescu, S.I., Huynh, J., Sim, F., and Thomas, M.
  (2020{\natexlab{a}}).
\newblock Partial pole placement via delay action: A python software for
  delayed feedback stabilizing design.
\newblock In \emph{2020 24th International Conference on System Theory, Control
  and Computing (ICSTCC)}, 196--201.

\bibitem[{Boussaada et~al.(2021)Boussaada, Mazanti, Niculescu, Leclerc, Raj,
  and Perraudin}]{Boussaada-TDS2021}
Boussaada, I., Mazanti, G., Niculescu, S.I., Leclerc, A., Raj, J., and
  Perraudin, M. (2021).
\newblock New features of {P}3delta software: Partial pole placement via delay
  action.
\newblock \emph{IFAC-PapersOnLine}, 54(18), 215--221.
\newblock 16th IFAC Workshop on Time Delay Systems.

\bibitem[{Boussaada and
  Niculescu(2016{\natexlab{a}})}]{Boussaada2016Characterizing}
Boussaada, I. and Niculescu, S.I. (2016{\natexlab{a}}).
\newblock Characterizing the codimension of zero singularities for time-delay
  systems: a link with {V}andermonde and {B}irkhoff incidence matrices.
\newblock \emph{Acta Appl. Math.}, 145, 47--88.

\bibitem[{Boussaada and Niculescu(2016{\natexlab{b}})}]{Boussaada2016Tracking}
Boussaada, I. and Niculescu, S.I. (2016{\natexlab{b}}).
\newblock Tracking the algebraic multiplicity of crossing imaginary roots for
  generic quasipolynomials: a {V}andermonde-based approach.
\newblock \emph{IEEE Trans. Automat. Control}, 61(6), 1601--1606.

\bibitem[{Boussaada and Niculescu(2018)}]{Boussaada2018Dominancy}
Boussaada, I. and Niculescu, S.I. (2018).
\newblock On the dominancy of multiple spectral values for time-delay systems
  with applications.
\newblock \emph{IFAC-PapersOnLine}, 51(14), 55--60.
\newblock 14th IFAC Workshop on Time Delay Systems.

\bibitem[{Boussaada et~al.(2020{\natexlab{b}})Boussaada, Niculescu, El-Ati,
  P\'{e}rez-Ramos, and Trabelsi}]{Boussaada2020Multiplicity}
Boussaada, I., Niculescu, S.I., El-Ati, A., P\'{e}rez-Ramos, R., and Trabelsi,
  K. (2020{\natexlab{b}}).
\newblock Multiplicity-induced-dominancy in parametric second-order delay
  differential equations: {A}nalysis and application in control design.
\newblock \emph{ESAIM Control Optim. Calc. Var.}, 26, Paper No. 57.

\bibitem[{Boussaada et~al.(2018{\natexlab{a}})Boussaada, Niculescu, and
  Trabelsi}]{Boussaada2018Towards}
Boussaada, I., Niculescu, S.I., and Trabelsi, K. (2018{\natexlab{a}}).
\newblock Towards a decay rate assignment based design for time-delay systems
  with multiple spectral values.
\newblock In \emph{Proceedings of the 23rd International Symposium on
  Mathematical Theory of Networks and Systems (MTNS)}, 864--871.

\bibitem[{Boussaada et~al.(2018{\natexlab{b}})Boussaada, Tliba, Niculescu,
  \"{U}nal, and Vyhl\'{\i}dal}]{Boussaada2018Further}
Boussaada, I., Tliba, S., Niculescu, S.I., \"{U}nal, H.U., and Vyhl\'{\i}dal,
  T. (2018{\natexlab{b}}).
\newblock Further remarks on the effect of multiple spectral values on the
  dynamics of time-delay systems. {A}pplication to the control of a mechanical
  system.
\newblock \emph{Linear Algebra Appl.}, 542, 589--604.

\bibitem[{Boussaada et~al.(2016)Boussaada, \"{U}nal, and
  Niculescu}]{Boussaada2016Multiplicity}
Boussaada, I., \"{U}nal, H.U., and Niculescu, S.I. (2016).
\newblock Multiplicity and stable varieties of time-delay systems: A missing
  link.
\newblock In \emph{Proceedings of the 22nd International Symposium on
  Mathematical Theory of Networks and Systems (MTNS)}, 188--194.

\bibitem[{Callender et~al.(1936)Callender, Hartree, and Porter}]{chp:36}
Callender, A., Hartree, D.R., and Porter, A. (1936).
\newblock Time-lag in a control system.
\newblock \emph{Phil. Trans. Royal Soc. Series A}, 235, 415--444.

\bibitem[{Gu et~al.(2003)Gu, Kharitonov, and Chen}]{Gu2003Stability}
Gu, K., Kharitonov, V.L., and Chen, J. (2003).
\newblock \emph{Stability of time-delay systems}.
\newblock Control Engineering. Birkh\"{a}user Boston, Inc., Boston, MA.

\bibitem[{Hale and Verduyn~Lunel(1993)}]{Hale1993Introduction}
Hale, J.K. and Verduyn~Lunel, S.M. (1993).
\newblock \emph{Introduction to functional differential equations}.
\newblock Springer-Verlag, New York.

\bibitem[{Hartree et~al.(1937)Hartree, Porter, Callender, and
  Stevenson}]{hpcs:37}
Hartree, D.R., Porter, A., Callender, A., and Stevenson, A. (1937).
\newblock Time-lag in a control system. {II}.
\newblock \emph{Proc. Royal Soc. A}, 161, 460--475.

\bibitem[{Ma et~al.(2022)Ma, Boussaada, Chen, Bonnet, Niculescu, and
  Chen}]{MBCBNC-Automatica-2022}
Ma, D., Boussaada, I., Chen, J., Bonnet, C., Niculescu, S.I., and Chen, J.
  (2022).
\newblock {PID} control design for first-order delay systems via {MID} pole
  placement: Performance vs. robustness.
\newblock \emph{Automatica}, 137, 110102.

\bibitem[{Mazanti et~al.(2020{\natexlab{a}})Mazanti, Boussaada, and
  Niculescu}]{Mazanti2020Qualitative}
Mazanti, G., Boussaada, I., and Niculescu, S.I. (2020{\natexlab{a}}).
\newblock On qualitative properties of single-delay linear retarded
  differential equations: Characteristic roots of maximal multiplicity are
  necessarily dominant.
\newblock \emph{IFAC-PapersOnLine}, 53(2), 4345--4350.
\newblock 21st IFAC WC.

\bibitem[{Mazanti et~al.(2021)Mazanti, Boussaada, and
  Niculescu}]{MazantiMultiplicity}
Mazanti, G., Boussaada, I., and Niculescu, S.I. (2021).
\newblock Multiplicity-induced-dominancy for delay-differential equations of
  retarded type.
\newblock \emph{J. Differential Equations}, 286, 84--118.

\bibitem[{Mazanti et~al.(2020{\natexlab{b}})Mazanti, Boussaada, Niculescu, and
  Vyhl{\'{\i}}dal}]{Mazanti2020Spectral}
Mazanti, G., Boussaada, I., Niculescu, S.I., and Vyhl{\'{\i}}dal, T.
  (2020{\natexlab{b}}).
\newblock Spectral dominance of complex roots for single-delay linear
  equations.
\newblock \emph{IFAC-PapersOnLine}, 53(2), 4357--4362.
\newblock 21st IFAC World Congress.

\bibitem[{Michiels and Niculescu(2014)}]{Michiels2014Stability}
Michiels, W. and Niculescu, S.I. (2014).
\newblock \emph{Stability, control, and computation for time-delay systems: An
  eigenvalue-based approach}.
\newblock SIAM, Philadelphia, PA, second edition.

\bibitem[{Molnar et~al.(2021)Molnar, Balogh, Boussaada, and
  Insperger}]{Molnar2021Calculation}
Molnar, C.A., Balogh, T., Boussaada, I., and Insperger, T. (2021).
\newblock Calculation of the critical delay for the double inverted pendulum.
\newblock \emph{J. Vib. Control}, 27(3-4), 356--364.

\bibitem[{Niculescu et~al.(2010)Niculescu, Michiels, Gu, and
  Abdallah}]{Niculescu2010Delay}
Niculescu, S.I., Michiels, W., Gu, K., and Abdallah, C.T. (2010).
\newblock Delay effects on output feedback control of dynamical systems.
\newblock In F.M. Atay (ed.), \emph{Complex time-delay systems}, 63--84.
  Springer, Berlin.

\bibitem[{Porter(1952)}]{porter:52}
Porter, A. (1952).
\newblock \emph{Introduction to servomechanisms}.
\newblock Wiley: London, 2nd edition.

\bibitem[{Sipahi et~al.(2011)Sipahi, Niculescu, Abdallah, Michiels, and
  Gu}]{Sipahi2011Stability}
Sipahi, R., Niculescu, S.I., Abdallah, C.T., Michiels, W., and Gu, K. (2011).
\newblock Stability and stabilization of systems with time delay: limitations
  and opportunities.
\newblock \emph{IEEE Control Syst. Mag.}, 31(1), 38--65.

\bibitem[{St\'{e}p\'{a}n(1989)}]{Stepan1989Retarded}
St\'{e}p\'{a}n, G. (1989).
\newblock \emph{Retarded dynamical systems: stability and characteristic
  functions}, volume 210 of \emph{Pitman Research Notes in Mathematics Series}.
\newblock Longman Scientific \& Technical, Harlow.

\bibitem[{Suh and Bien(1979)}]{Suh1979Proportional}
Suh, I.H. and Bien, Z. (1979).
\newblock Proportional minus delay controller.
\newblock \emph{IEEE Trans. Automat. Control}, 24(2), 370--372.

\end{thebibliography}

\end{document}